# FUNCTIONAL LINEAR REGRESSION THAT'S INTERPRETABLE[1]

By Gareth M. James, Jing Wang and Ji Zhu

*University of Southern California, University of Michigan
and University of Michigan*

Regression models to relate a scalar $Y$ to a functional predictor $X(t)$ are becoming increasingly common. Work in this area has concentrated on estimating a coefficient function, $\beta(t)$, with $Y$ related to $X(t)$ through $\int \beta(t)X(t)\,dt$. Regions where $\beta(t) \neq 0$ correspond to places where there is a relationship between $X(t)$ and $Y$. Alternatively, points where $\beta(t) = 0$ indicate no relationship. Hence, for interpretation purposes, it is desirable for a regression procedure to be capable of producing estimates of $\beta(t)$ that are exactly zero over regions with no apparent relationship and have simple structures over the remaining regions. Unfortunately, most fitting procedures result in an estimate for $\beta(t)$ that is rarely exactly zero and has unnatural wiggles making the curve hard to interpret. In this article we introduce a new approach which uses variable selection ideas, applied to various derivatives of $\beta(t)$, to produce estimates that are both interpretable, flexible and accurate. We call our method "Functional Linear Regression That's Interpretable" (FLiRTI) and demonstrate it on simulated and real-world data sets. In addition, non-asymptotic theoretical bounds on the estimation error are presented. The bounds provide strong theoretical motivation for our approach.

**1. Introduction.** In recent years functional data analysis (FDA) has become an increasingly important analytical tool as more data has arisen where the primary unit of observation can be viewed as a curve or in general a function. One of the most useful tools in FDA is that of functional regression. This setting can correspond to either functional predictors or functional responses. See Ramsay and Silverman (2002) and Muller and Stadtmuller (2005) for numerous specific applications. One commonly studied problem

Received April 2008; revised July 2008.

[1]Supported by NSF Grants DMS-07-05312, DMS-05-05432 and DMS-07-05532.

involves data containing functional responses. A sampling of papers examining this situation includes Fahrmeir and Tutz (1994), Liang and Zeger (1986), Faraway (1997), Hoover et al. (1998), Wu et al. (1998), Fan and Zhang (2000) and Lin and Ying (2001). However, in this paper, we are primarily interested in the alternative situation, where we obtain a set of observations $\{X_i(t), Y_i\}$ for $i = 1, \ldots, n$, where $X_i(t)$ is a functional predictor and $Y_i$ a real valued response. Ramsay and Silverman (2005) discuss this scenario and several papers have also been written on the topic, both for continuous and categorical responses, and for linear and nonlinear models [Hastie and Mallows (1993), James and Hastie (2001), Ferraty and Vieu (2002), James (2002), Ferraty and Vieu (2003), Muller and Stadtmuller (2005), James and Silverman (2005)].

Since our primary interest here is interpretation, we will be examining the standard functional linear regression (FLR) model, which relates functional predictors to a scalar response via

$$(1) \qquad Y_i = \beta_0 + \int X_i(t)\beta(t)\,dt + \varepsilon_i, \qquad i = 1, \ldots, n,$$

where $\beta(t)$ is the "coefficient function." We will assume that $X_i(t)$ is scaled such that $0 \leq t \leq 1$. Clearly, for any finite $n$, it would be possible to perfectly interpolate the responses if no restrictions were placed on $\beta(t)$. Such restrictions generally take one of two possible forms. The first method, which we call the "basis approach," involves representing $\beta(t)$ using a $p$-dimensional basis function, $\beta(t) = \mathbf{B}(t)^T \boldsymbol{\eta}$, where $p$ is hopefully large enough to capture the patterns in $\beta(t)$ but small enough to regularize the fit. With this method (1) can be reexpressed as $Y_i = \beta_0 + \mathbf{X}_i^T \boldsymbol{\eta} + \varepsilon_i$, where $\mathbf{X}_i = \int X_i(t)\mathbf{B}(t)\,dt$, and $\boldsymbol{\eta}$ can be estimated using ordinary least squares. The second method, which we call the "penalization approach," involves a penalized least squares estimation procedure to shrink variability in $\beta(t)$. A standard penalty is of the form $\int \beta^{(d)}(t)^2\,dt$ with $d = 2$ being a common choice. In this case one would find $\beta(t)$ to minimize $\sum_{i=1}^{n}(Y_i - \beta_0 - \int X_i(t)\beta(t)\,dt)^2 + \lambda \int \beta^{(d)}(t)^2\,dt$ for some $\lambda > 0$.

As with standard linear regression, $\beta(t)$ determines the effect of $X_i(t)$ on $Y_i$. For example, changes in $X_i(t)$ have no effect on $Y_i$ over regions, where $\beta(t) = 0$. Alternatively, changes in $X_i(t)$ have a greater effect on $Y_i$ over regions, where $|\beta(t)|$ is large. Hence, in terms of interpretation, coefficient curves with certain structures are more appealing than others. For example, if $\beta(t)$ is exactly zero over large regions then $X_i(t)$ only has an effect on $Y_i$ over the remaining time points. Additionally, if $\beta(t)$ is constant for any given non-zero region then the effect of $X_i(t)$ on $Y_i$ remains constant within that region. Finally, if $\beta(t)$ is exactly linear for any given region then the change in the effect of $X_i(t)$ is constant over that region. Clearly the interpretation of the predictor-response relationship is more difficult as the shape of $\beta(t)$ becomes more complicated. Unfortunately, the basis and penalization approaches both generate $\beta(t)$ curves that exhibit wiggles and are not exactly linear or constant over any region. In addition, $\beta(t)$ will be exactly equal to zero at no more than a few locations even if there is no relationship between $X(t)$ and $Y$ for large regions of $t$.

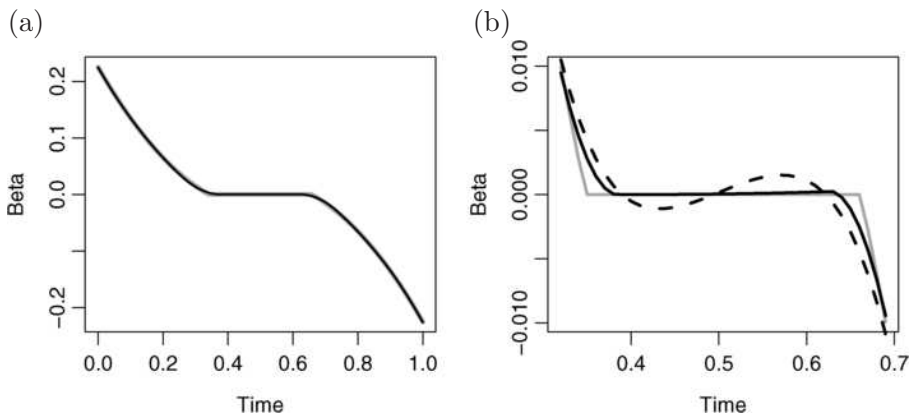

Fig. 1. (a) *True beta curve (grey) generated from two quadratic curves and a section with $\beta(t) = 0$. The FLiRTI estimate from constraining the zeroth and third derivative is shown in black.* (b) *Same plot for the region $0.3 \leq t \leq 0.7$. The dashed line is the best B-spline fit.*

In this paper we develop a new method, which we call "Functional Linear Regression That's Interpretable" (FLiRTI), that produces accurate but also highly interpretable estimates for the coefficient function $\beta(t)$. Additionally, it is computationally efficient, extremely flexible in terms of the form of the estimate and has highly desirable theoretical properties. The key to our procedure is to reformulate the problem as a form of variable selection. In particular we divide the time period up into a fine grid of points. We then use variable selection methods to determine whether the $d$th derivative of $\beta(t)$ is zero or not at each of the grid points. The implicit assumption is that the $d$th derivative will be zero at most grid points, that is, it will be sparse. By choosing appropriate derivatives one can produce a large range of highly interpretable $\beta(t)$ curves. Consider, for example, Figures 1–3, which illustrate a range of FLiRTI fits applied to simulated data sets. In Figure 1 the true $\beta(t)$ used to generate the data consisted of a quadratic curve and a section with $\beta(t) = 0$. Figure 1(a) plots the FLiRTI estimate, produced by assuming sparsity in the zeroth and third derivatives. The sparsity in the zeroth derivative generates the zero section while the sparsity in the third derivative ensures a smooth fit. Figure 1(b) illustrates the same plot concentrating on the region between 0.3 and 0.7. Notice that the corresponding B-spline estimate, represented by the dashed line, provides a poor approximation for the region, where $\beta(t) = 0$. It is important to note that we did not specify which regions would have zero derivatives, the FLiRTI procedure is capable of automatically selecting the appropriate shape. In Figure 2 $\beta(t)$ was chosen as a piecewise linear curve with the middle section set to zero. Figure 2(a) shows the corresponding FLiRTI estimate generated by assuming sparsity in the zeroth and second derivatives and is almost a perfect fit. Alternatively, one can produce a smoother, but slightly less easily interpretable fit, by assuming sparsity in higher-order derivatives. Figure 2(b), which concentrates on the region between $t = 0.2$ to $t = 0.8$, plots the FLiRTI fit assuming sparsity in the zeroth and third derivative. Notice that the sparsity in the third derivative induces a smoother estimate with little sacrifice in accuracy. Finally, Figure 3 illustrates a FLiRTI fit applied to data generated using a simple cubic $\beta(t)$ curve, a situation where one might not expect FLiRTI to provide any advantage over a standard approach such as using a B-spline basis. However, the figure, along with the simulation results in Section 5, shows that even in this situation the FLiRTI method gives highly accurate estimates. These three examples illustrate the flexibility of FLiRTI in that it can produce estimates ranging from highly interpretable simple linear fits, through smooth fits with zero regions, to more complicated nonlinear structures with equal ease. The key idea here is that, given a strong signal in the data, FLiRTI is flexible enough to estimate $\beta(t)$ accurately. However, in situations where the signal is weaker, FLiRTI will automatically shrink the estimated $\beta(t)$ towards a more interpretable structure.

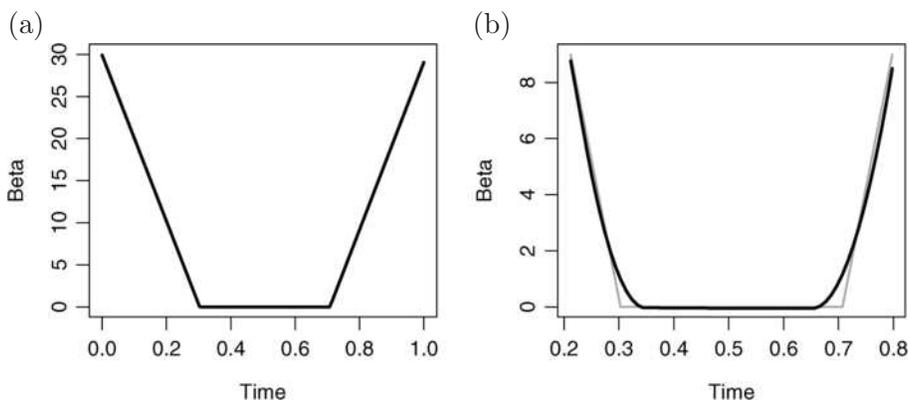

Fig. 2. *Plots of true beta curve (grey) and corresponding FLiRTI estimates (black). For each plot we constrained* (a) *zeroth and second derivative,* (b) *zeroth and third derivative.*

The paper is laid out as follows. In Section 2 we develop the FLiRTI model and also detail two fitting procedures, one making use of the lasso [Tibshirani (1996)] and the other utilizing the Dantzig selector [Candes and Tao (2007)]. The theoretical developments for our method are presented in Section 3

(a) (b)

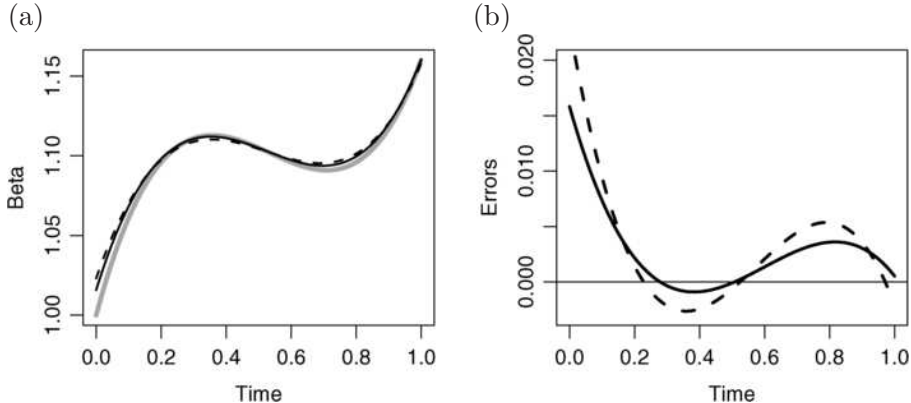

FIG. 3. (a) *True beta curve (grey) generated from a cubic curve. The FLiRTI estimate from constraining the zeroth and fourth derivative is represented by the solid black line and the B-spline estimate is the dashed line.* (b) *Estimation errors using the B-spline method (dashed) and FLiRTI (solid).*

where we outline both nonasymptotic bounds on the error as well as asymptotic properties of our estimate as $n$ grows. Then, in Section 4, we extend the FLiRTI method in two directions. First we show how to control multiple derivatives simultaneously, which allows us to, for example, produce a $\beta(t)$ curve that is exactly zero in certain sections and exactly linear in other sections. Second, we develop a version of FLiRTI that can be applied to generalized linear models (GLM). A detailed simulation study is presented in Section 5. Finally, we apply FLiRTI to real world data in Section 6 and end with a discussion in Section 7.

**2. FLiRTI methodology.** In this section we first develop the FLiRTI model and then demonstrate how we use the lasso and Dantzig selector methods to fit it and hence estimate $\beta(t)$.

2.1. *The FLiRTI model.* Our approach borrows ideas from the basis and penalization methods but is rather different from either. We start in a similar vein to the basis approach by selecting a $p$-dimensional basis $\mathbf{B}(t) = [b_1(t), b_2(t), \ldots, b_p(t)]^T$. However, instead of assuming $\mathbf{B}(t)$ provides a perfect fit for $\beta(t)$, we allow for some error using the model

$$(2) \qquad \beta(t) = \mathbf{B}(t)^T \boldsymbol{\eta} + e(t),$$

where $e(t)$ represents the deviations of the true $\beta(t)$ from our model. In addition, unlike the basis approach where $p$ is chosen small to provide some form of regularization, we typically choose $p \gg n$ so $|e(t)|$ can generally be assumed to be small. In Section 3 we show that the error in the estimate for $\beta(t)$, that is, $|\hat{\beta}(t) - \beta(t)|$, can potentially be of order $\sqrt{\log(p)/n}$. Hence, low error rates can be achieved even for values of $p$ much larger than $n$. Our theoretical results apply to any high dimensional basis, such as splines, Fourier or wavelets. For the empirical results presented in this paper we opted to use a simple grid basis, where $b_k(t)$ equals 1 if $t \in R_k = \{t : \frac{k-1}{p} < t \leq \frac{k}{p}\}$ and 0 otherwise.

Combining (1) and (2) we arrive at

$$(3) \qquad\qquad Y_i = \beta_0 + \mathbf{X}_i^T \boldsymbol{\eta} + \varepsilon_i^*,$$

where $\mathbf{X}_i = \int X_i(t)\mathbf{B}(t)\,dt$ and $\varepsilon_i^* = \varepsilon_i + \int X_i(t)e(t)\,dt$. Estimating $\boldsymbol{\eta}$ presents a difficulty because $p > n$. One could potentially estimate $\boldsymbol{\eta}$ using a variable selection procedure except that for an arbitrary basis, $\mathbf{B}(t)$, there is no reason to suppose that $\boldsymbol{\eta}$ will be sparse. In fact for many bases $\boldsymbol{\eta}$ will contain no zero elements. Instead we model $\beta(t)$ assuming that one or more of its derivatives are sparse, that is, $\beta^{(d)}(t) = 0$ over large regions of $t$ for one or more values of $d = 0, 1, 2, \ldots$. This model has the advantage of both constraining $\boldsymbol{\eta}$ enough to allow us to fit (3) as well as producing a highly interpretable estimate for $\beta(t)$. For example, $\beta^{(0)}(t) = 0$ guarantees $X(t)$ has no effect on $Y$ at $t$, $\beta^{(1)}(t) = 0$ implies that $\beta(t)$ is constant at $t$, $\beta^{(2)}(t) = 0$ means that $\beta(t)$ is linear at $t$, etc.

Let $A = [D^d\mathbf{B}(t_1), D^d\mathbf{B}(t_2), \ldots, D^d\mathbf{B}(t_p)]^T$, where $t_1, t_2, \ldots, t_p$ represent a grid of $p$ evenly spaced points and $D^d$ is the $d$th finite difference operator, that is, $D\mathbf{B}(t_j) = p[\mathbf{B}(t_j) - \mathbf{B}(t_{j-1})]$, $D^2\mathbf{B}(t_j) = p^2[\mathbf{B}(t_j) - 2\mathbf{B}(t_{j-1}) + \mathbf{B}(t_{j-2})]$, etc. Then, if

$$(4) \qquad\qquad \boldsymbol{\gamma} = A\boldsymbol{\eta},$$

$\gamma_j$ provides an approximation to $\beta^{(d)}(t_j)$ and hence, enforcing sparsity in $\gamma$ constrains $\beta^{(d)}(t_j)$ to be zero at most time points. For example, one may believe that $\beta^{(2)}(t) = 0$ over many regions of $t$, that is, $\beta(t)$ is exactly linear over large regions of $t$. In this situation we would let

$$(5) \qquad\qquad A = [D^2\mathbf{B}(t_1), D^2\mathbf{B}(t_2), \ldots, D^2\mathbf{B}(t_p)]^T,$$

which implies $\gamma_j = p^2[\mathbf{B}(t_j)^T\boldsymbol{\eta} - 2\mathbf{B}(t_{j-1})^T\boldsymbol{\eta} + \mathbf{B}(t_{j-2})^T\boldsymbol{\eta}]$. Hence, provided $p$ is large, so $t$ is sampled on a fine grid, and $e(t)$ is smooth, $\gamma_j \approx \beta^{(2)}(t_j)$. In this case enforcing sparsity in the $\gamma_j$'s will produce an estimate for $\beta(t)$ that is linear except at the time points corresponding to nonzero values of $\gamma_j$.

If $A$ is constructed using a single derivative, as in (5), then we can always choose a grid of $p$ different time points, $t_1, t_2, \ldots, t_p$ such that $A$ is a square $p$ by $p$ invertible matrix. In this case $\boldsymbol{\eta} = A^{-1}\boldsymbol{\gamma}$ so we may combine (3) and (4) to produce the FLiRTI model

$$(6) \qquad\qquad \mathbf{Y} = V\boldsymbol{\gamma} + \boldsymbol{\varepsilon}^*,$$

where $V = [\mathbf{1}|XA^{-1}]$, $\mathbf{1}$ is a vector of ones and $\beta_0$ has been incorporated into $\boldsymbol{\gamma}$. We discuss the situation with multiple derivatives, where $A$ may no longer be invertible, in Section 4.

2.2. *Fitting the model.* Since $\boldsymbol{\gamma}$ is assumed sparse, one could potentially use a variety of variable selection methods to fit (6). There has recently been a great deal of development of new model selection methods that work with large values of $p$. A few examples include the lasso [Tibshirani (1996), Chen, Donoho and Saunders (1998)], SCAD [Fan and Li (2001)], the Elastic Net [Zou and Hastie (2005)], the Dantzig selector [Candes and Tao (2007)] and VISA [Radchenko and James (2008)]. We opted to explore both the lasso and Dantzig selector for several reasons. First, both methods have demonstrated strong empirical results on models with large values of $p$. Second, the LARS algorithm [Efron et al. (2004)] can be used to efficiently compute the coefficient path for the lasso. Similarly the DASSO algorithm [James, Radchenko and Lv (2009)], which is a generalization of LARS, will efficiently compute the Dantzig selector coefficient path. Finally, we demonstrate in Section 3 that identical non-asymptotic bounds can be placed on the errors in the estimates of $\beta(t)$ that result from either approach.

Consider the linear regression model $\mathbf{Y} = X\boldsymbol{\beta} + \boldsymbol{\varepsilon}$. Then the lasso estimate, $\widehat{\boldsymbol{\beta}}_{\mathrm{L}}$, is defined by

$$(7) \qquad \widehat{\boldsymbol{\beta}}_{\mathrm{L}} = \arg\min_{\boldsymbol{\beta}} \frac{1}{2}\|\mathbf{Y} - X\boldsymbol{\beta}\|_2^2 + \lambda\|\boldsymbol{\beta}\|_1,$$

where $\|\cdot\|_1$ and $\|\cdot\|_2$ respectively denote the $L_1$ and $L_2$ norms and $\lambda \geq 0$ is a tuning parameter. Alternatively, the Dantzig selector estimate, $\widehat{\boldsymbol{\beta}}_{\mathrm{DS}}$, is given by

$$(8) \qquad \widehat{\boldsymbol{\beta}}_{\mathrm{DS}} = \arg\min_{\boldsymbol{\beta}} \|\boldsymbol{\beta}\|_1 \quad \text{subject to} \quad |\mathbf{X}_j^T(\mathbf{Y} - X\boldsymbol{\beta})| \leq \lambda, \quad j = 1, \ldots, p,$$

where $\mathbf{X}_j$ is the $j$th column of $X$ and $\lambda \geq 0$ is a tuning parameter.

Using either approach the LARS or DASSO algorithms can be used to efficiently compute all solutions for various values of the tuning parameter, $\lambda$. Hence, using a validation/cross-validation approach to select $\lambda$ can be easily implemented. To generate the final FLiRTI estimate we first produce $\widehat{\boldsymbol{\gamma}}$ by fitting the FLiRTI model, (6), using either the lasso, (7) or the Dantzig selector, (8). Note that both methods generally assume a standardized design matrix with columns of norm one. Hence, we first standardize $V$, apply the lasso or Dantzig selector and then divide the resulting coefficient estimates by the original column norms of $V$ to produce $\widehat{\boldsymbol{\gamma}}$. After the coefficients, $\widehat{\boldsymbol{\gamma}}$, have been obtained we produce the FLiRTI estimate for $\beta(t)$ using

$$(9) \qquad \widehat{\beta}(t) = \mathbf{B}(t)^T\widehat{\boldsymbol{\eta}} = \mathbf{B}(t)^T A^{-1}\widehat{\boldsymbol{\gamma}}_{(-1)},$$

where $\widehat{\boldsymbol{\gamma}}_{(-1)}$ represents $\widehat{\boldsymbol{\gamma}}$ after removing the estimate for $\beta_0$.

**3. Theoretical results.** In this section we show that not only does the FLiRTI approach empirically produce good estimates for $\beta(t)$, but that for any $p$ by $p$ invertible $A$ we can in fact prove tight, nonasymptotic bounds on the error in our estimate. In addition, we derive asymptotic rates of convergence. Note that for notational convenience the results in this section assume $\beta_0 = 0$ and drop the intercept term from the model. However, the theory all extends in a straightforward manner to the situation with $\beta_0$ unknown.

3.1. *A nonasymptotic bound on the error.* Let $\widehat{\boldsymbol{\gamma}}_\lambda$ correspond to the lasso solution using tuning parameter $\lambda$. Let $D_\lambda$ be a diagonal matrix with $j$th diagonal equal to $1, -1$ or $0$ depending on whether the $j$th component of $\widehat{\boldsymbol{\gamma}}_\lambda$ is positive, negative or zero, respectively. Consider the following condition on the design matrix, $V$.

$$(10) \qquad \mathbf{u} = (D_\lambda \tilde{V}^T \tilde{V} D_\lambda)^{-1} \mathbf{1} \geq \mathbf{0} \quad \text{and} \quad \|\tilde{V}^T \tilde{V} D_\lambda \mathbf{u}\|_\infty \leq 1,$$

where $\tilde{V}$ corresponds to $V$ after standardizing its columns, $\mathbf{1}$ is a vector of ones and the inequality for vectors is understood componentwise. James, Radchenko and Lv (2009) prove that when (10) holds the Dantzig selector's nonasymptotic bounds [Candes and Tao (2007)] also apply for the lasso. Our Theorem 1 makes use of (10) to provide a nonasymptotic bound on the $L_2$ error in the FLiRTI estimate, using either the Dantzig selector or the lasso. Note that the values $\delta, \theta$ and $C_{n,p}(t)$ are all known constants which we have defined in the proof of this result provided in Appendix A.

THEOREM 1. *For a given $p$-dimensional basis $\mathbf{B}_p(t)$, let $\omega_p = \sup_t |e_p(t)|$ and $\boldsymbol{\gamma}_p = A\boldsymbol{\eta}_p$, where $A$ is a $p$ by $p$ matrix. Suppose that $\boldsymbol{\gamma}_p$ has at most $S_p$ non-zero components and $\delta_{2S_p}^V + \theta_{S_p, 2S_p}^V < 1$. Further, suppose that we estimate $\beta(t)$ using the FLiRTI estimate given by (9) using any value of $\lambda$ such that (10) holds and*

$$(11) \qquad \max |\tilde{V}^T \boldsymbol{\varepsilon}^*| \leq \lambda.$$

*Then, for every $0 \leq t \leq 1$,*

$$(12) \qquad |\widehat{\beta}(t) - \beta(t)| \leq \frac{1}{\sqrt{n}} C_{n,p}(t) \lambda \sqrt{S_p} + \omega_p.$$

The constants $\delta$ and $\theta$ are both measures of the orthogonality of $V$. The closer they are to zero the closer $V$ is to orthogonal. The condition $\delta_{2S_p}^V + \theta_{S_p, 2S_p}^V < 1$, which is utilized in the paper of Candes and Tao (2007), ensures that $\beta(t)$ is identifiable. It should be noted that (10) is only required when using the lasso to compute FLiRTI. The above results hold for the Dantzig selector even if (10) is violated. While this is a slight theoretical advantage for the Dantzig selector, our simulation results suggest that both methods perform well in practice. Theorem 1 suggests that our optimal choice for $\lambda$ would be the lowest value such that (11) holds. Theorem 2 shows how to choose $\lambda$ such that (11) holds with high probability.

THEOREM 2. *Suppose that* $\varepsilon_i \sim N(0, \sigma_1^2)$ *and that there exits an* $M < \infty$ *such that* $\int |X_i(t)| \, dt \leq M$ *for all* $i$. *Then for any* $\phi \geq 0$, *if* $\lambda = \sigma_1 \sqrt{2(1+\phi) \log p} + M \omega_p \sqrt{n}$ *then* (11) *will hold with probability at least* $1 - (p^\phi \sqrt{4\pi(1+\phi) \log p})^{-1}$, *and hence*

$$
\begin{aligned}
(13) \qquad |\widehat{\beta}(t) - \beta(t)| \leq \frac{1}{\sqrt{n}} &C_{n,p}(t) \sigma_1 \sqrt{2 S_p (1+\phi) \log p} \\
&+ \omega_p \{1 + C_{n,p}(t) M \sqrt{S_p}\}.
\end{aligned}
$$

*In addition, if we assume* $\varepsilon_i^* \sim N(0, \sigma_2^2)$ *then* (11) *will hold with the same probability for* $\lambda = \sigma_2 \sqrt{2(1+\phi) \log p}$ *in which case*

$$
(14) \qquad |\widehat{\beta}(t) - \beta(t)| \leq \frac{1}{\sqrt{n}} C_{n,p}(t) \sigma_2 \sqrt{2 S_p (1+\phi) \log p} + \omega_p.
$$

Note that (13) and (14) are non-asymptotic results that hold, with high probability, for any $n$ or $p$. One can show that, under suitable conditions, $C_{n,p}(t)$ converges to a constant as $n$ and $p$ grow. In this case the first terms of (13) and (14) are proportional to $\sqrt{\log p / n}$ while the second terms will generally shrink as $\omega_p$ declines with $p$. For example, using the piecewise constant basis it is easy to show that $\omega_p$ converges to zero at a rate of $1/p$ provided $\beta'(t)$ is bounded. Alternatively using a piecewise polynomial basis of order $d$ then $\omega_p$ converges to zero at a rate of $1/p^{d+1}$ provided $\beta^{(d+1)}(t)$ is bounded.

3.2. *Asymptotic rates of convergence.* The bounds presented in Theorem 1 can be used to derive asymptotic rates of convergence for $\widehat{\beta}(t)$ as $n$ and $p$ grow. The exact convergence rates are dependent on the choice of $\mathbf{B}_p(t)$ and $A$ so we first state A-1 through A-6, which give general conditions for convergence. We show in Theorems 3–6 that these conditions are sufficient to guarantee convergence for any choice of $\mathbf{B}_p(t)$ and $A$, and then Corollary 1 provides specific examples where the conditions can be shown to hold. Let $\alpha_{n,p}(t) = (1 - \delta_{2S}^{V_{n,p}} + \theta_{S,2S}^{V_{n,p}}) C_{n,p}$.

A-1 There exists $S < \infty$ such that $S_p \leq S$ for all $p$.

A-2 There exists $m > 0$ such that $p^m \omega_p$ is bounded, that is, $\omega_p \leq H/p^m$ for some $H < \infty$.

A-3 For a given $t$, there exists $b_t$ such that $p^{-b_t}\alpha_{n,p}(t)$ is bounded for all $n$ and $p$.

A-4 There exists $c$ such that $p^{-c}\sup_t \alpha_{n,p}(t)$ is bounded for all $n$ and $p$.

A-5 There exists a $p^*$ such that $\delta_{2S}^{V_{n,p^*}} + \theta_{S,2S}^{V_{n,p^*}}$ is bounded away from one for large enough $n$.

A-6 $\delta_{2S}^{V_{n,p_n}} + \theta_{S,2S}^{V_{n,p_n}}$ is bounded away from one for large enough $n$, where $n \to \infty, p_n \to \infty$ and $p_n/n \to 0$.

A-1 states that the number of changes in the derivative of $\beta(t)$ is bounded. A-2 assumes that the bias in our estimate for $\beta(t)$ converges to zero at the rate of $p^m$, for some $m > 0$. A-3 requires that $\alpha_{n,p_n}(t)$ grows no faster than $p^{b_t}$. A-4 is simply a stronger form of A-3. A-5 and A-6 both ensure that the design matrix is close enough to orthogonal for (13) to hold and hence imposes a form of identifiability on $\beta(t)$. For the following two theorems we assume that the conditions in Theorem 1 hold, $\lambda$ is set to $\sigma_1\sqrt{2(1+\phi)\log p} + M\omega_p\sqrt{n}$ and $\varepsilon_i \sim N(0, \sigma_1^2)$.

THEOREM 3. *Suppose* A-1 *through* A-5 *all hold and we fix* $p = p^*$. *Then, with arbitrarily high probability, as* $n \to \infty$,

$$|\widehat{\beta}_n(t) - \beta(t)| \le O(n^{-1/2}) + E_n(t)$$

*and*

$$\sup_t |\widehat{\beta}_n(t) - \beta(t)| \le O(n^{-1/2}) + \sup_t E_n(t),$$

*where* $E_n(t) = \frac{H}{p^{*m}}\{1 + C_{n,p^*}(t)M\sqrt{S}\}$.

More specifically, the probability referred to in Theorem 3 converges to one as $\phi \to \infty$. Theorem 3 states that, with our weakest set of assumptions, fixing $p$ as $n \to \infty$ will cause the FLiRTI estimate to be asymptotically within $E_n(t)$ of the true $\beta(t)$. $E_n(t)$ represents the bias in the approximation caused by representing $\beta(t)$ using a $p$ dimensional basis.

THEOREM 4. *Suppose we replace* A-5 *with* A-6. *Then, provided* $b_t$ *and* $c$ *are less than* $m$, *if we let* $p$ *grow at the rate of* $n^{1/(2m)}$,

$$|\widehat{\beta}_n(t) - \beta(t)| = O\left(\frac{\sqrt{\log n}}{n^{1/2 - b_t/(2m)}}\right)$$

*and*

$$\sup_t |\widehat{\beta}(t) - \beta(t)| = O\left(\frac{\sqrt{\log n}}{n^{1/2 - c/(2m)}}\right).$$

Theorem 4 shows that, assuming A-6 holds, $\widehat{\beta}_n(t)$ will converge to $\beta(t)$ at the given convergence rate. With additional assumptions, stronger results are possible. In the Appendix we present Theorems 5 and 6, which provide faster rates of convergence under the additional assumption that $\varepsilon_i^*$ has a mean zero Gaussian distribution.

Theorems 3–6 make use of assumptions A-1 to A-6. Whether these assumptions hold in practice depends on the choice of basis function and $A$ matrix. Corollary 1 below provides one specific example where conditions A-1 to A-4 can be shown to hold.

COROLLARY 1. *Suppose we divide the time interval* $[0,1]$ *into $p$ equal regions and use the piecewise constant basis. Let $A$ be the second difference matrix given by (23). Suppose that $X_i(t)$ is bounded above zero for all $i$ and $t$. Then, provided $\beta'(t)$ is bounded and $\beta''(t) \neq 0$ at a finite number of points,* A-1, A-2 *and* A-3 *all hold with $m = 1$, $b_0 = 0$ and $b_t = 0.5$, $0 < t < 1$. In addition, for $t$ bounded away from one,* A-4 *will also hold with $c = 0.5$. Hence, if* A-5 *holds and $\varepsilon_i^* \sim N(0, \sigma_2^2)$,*

$$|\widehat{\beta}_n(t) - \beta(t)| \leq O(n^{-1/2}) + \frac{H}{p^*} \quad and \quad \sup_t |\widehat{\beta}_n(t) - \beta(t)| \leq O(n^{-1/2}) + \frac{H}{p^*}.$$

*Alternatively, if* A-6 *holds and $\varepsilon_i^* \sim N(0, \sigma_2^2)$,*

$$|\widehat{\beta}_n(t) - \beta(t)| = \begin{cases} O\left(\dfrac{\sqrt{\log n}}{n^{1/2}}\right), & t = 0, \\ O\left(\dfrac{\sqrt{\log n}}{n^{1/3}}\right), & 0 < t < 1, \end{cases}$$

*and*

$$\sup_{0 < t < 1-a} |\widehat{\beta}_n(t) - \beta(t)| = O\left(\frac{\sqrt{\log n}}{n^{1/3}}\right)$$

*for any $a > 0$. Similar, though slightly weaker, results hold when $\varepsilon_i^*$ does not have a Gaussian distribution.*

Note that the choice of $t = 0$ for the faster rate of convergence is simply made for notational convenience. By appropriately choosing A we can achieve this rate for any fixed value of $t$ or indeed for any finite set of time points. In addition, the choice of the piecewise constant basis was made for simplicity. Similar results can be derived for higher-order polynomial bases in which case A-2 will hold with a higher $m$ and hence faster rates of convergence will be possible.

**4. Extensions.** In this section we discuss two useful extensions of the basic FLiRTI methodology. First, in Section 4.1, we show how to control multiple derivatives simultaneously to allow more flexibility in the types of possible shapes one can produce. Second, in Section 4.2, we extend FLiRTI to GLM models.

4.1. *Controlling multiple derivatives.* So far we have concentrated on controlling a single derivative of $\beta(t)$. However, one of the most powerful aspects of the FLiRTI approach is that we can combine constraints for multiple derivatives together to produce curves with many different properties. For example, one may believe that both $\beta^{(0)}(t) = 0$ and $\beta^{(2)}(t) = 0$ over many regions of $t$, that is, $\beta(t)$ is exactly zero over certain regions and $\beta(t)$ is exactly linear over other regions of $t$. In this situation, we would let

$$A = [D^0 \mathbf{B}(t_1), D^0 \mathbf{B}(t_2), \ldots, D^0 \mathbf{B}(t_p), D^2 \mathbf{B}(t_1), D^2 \mathbf{B}(t_2), \ldots, D^2 \mathbf{B}(t_p)]^T.$$
(15)

In general, such a matrix will have more rows than columns so will not be invertible. Let $A_{(1)}$ represent the first $p$ rows of $A$ and $A_{(2)}$ the remainder. Similarly, let $\boldsymbol{\gamma}_{(1)}$ represent the first $p$ elements of $\boldsymbol{\gamma}$ and $\boldsymbol{\gamma}_{(2)}$ the remaining elements. Then, assuming $A$ is arranged so that $A_{(1)}$ is invertible, the constraint $\boldsymbol{\gamma} = A\boldsymbol{\eta}$ implies

$$(16) \qquad \boldsymbol{\eta} = A_{(1)}^{-1}\boldsymbol{\gamma}_{(1)} \quad \text{and} \quad \boldsymbol{\gamma}_{(2)} = A_{(2)}A_{(1)}^{-1}\boldsymbol{\gamma}_{(1)}.$$

Hence, (3) can be expressed as

$$(17) \qquad Y_i = \beta_0 + (A_{(1)}^{-1T}\mathbf{X}_i)^T\boldsymbol{\gamma}_{(1)} + \varepsilon_i^*, \qquad i = 1, \ldots, n.$$

We then use this model to estimate $\boldsymbol{\gamma}$ subject to the constraint given by (16). We achieve this by implementing the Dantzig selector or lasso in a similar fashion to that in Section 2 except that we replace the old design matrix with $V_{(1)} = [\mathbf{1}|XA_{(1)}^{-1}]$ and (16) is enforced in addition to the usual constraints.

Finally, when constraining multiple derivatives one may well not wish to place equal weight on each derivative. For example, for the $A$ given by (15), we may wish to place a greater emphasis on sparsity in the second derivative than in the zeroth, or vice versa. Hence, instead of simply minimizing the $L_1$ norm of $\boldsymbol{\gamma}$ we minimize $\|\Omega\boldsymbol{\gamma}\|_1$, where $\Omega$ is a diagonal weighting matrix. In theory a different weight could be chosen for each $\gamma_j$ but in practice this would not be feasible. Instead, for an $A$ such as (15), we place a weight of one on the second derivatives and select a single weight, $\omega$, chosen via cross-validation, for the zeroth derivatives. This approach provides flexibility while still being computationally feasible and has worked well on all the problems we have examined. The FLiRTI fits in Figures 1–3 were produced using this multiple derivative methodology.

4.2. *FLiRTI for functional generalized linear models.* The FLiRTI model can easily be adapted to GLM data where the response is no longer assumed to be Gaussian. James and Radchenko (2009) demonstrate that the Dantzig selector can be naturally extended to the GLM domain by optimizing

$$(18) \qquad \min_{\boldsymbol{\beta}} \|\boldsymbol{\beta}\|_1 \quad \text{subject to} \quad |\mathbf{X}_j^T(\mathbf{Y} - \boldsymbol{\mu})| \le \lambda, \qquad j = 1, \dots, p,$$

where $\boldsymbol{\mu} = g^{-1}(X\boldsymbol{\beta})$ and $g$ is the canonical link function. Optimizing (18) no longer involves linear constraints but it can be solved using an iterative weighted linear programming approach. In particular, let $U_i$ be the conditional variance of $Y_i$ given the current estimate $\hat{\beta}$ and let $Z_i = \sum_{j=1}^p X_{ij}\hat{\beta} + (Y_i - \hat{\mu}_i)/U_i$. Then one can apply the standard Dantzig selector, using $Z_i^* = Z_i\sqrt{U_i}$ as the response and $X_{ij}^* = X_{ij}\sqrt{U_i}$ as the predictor, to produce a new estimate for $\hat{\beta}$. James and Radchenko (2009) show that iteratively applying the Dantzig selector in this fashion until convergence generally does a good job of solving (18). We apply the same approach, except that we iteratively apply the modified version of the Dantzig selector, outlined in Section 4.1, to the transformed response and predictor variables, using $V_{(1)} = [\mathbf{1}|X A_{(1)}^{-1}]$ as the design matrix. James and Radchenko (2009) also suggest an algorithm for approximating the GLM Dantzig selector coefficient path for different values of $\lambda$. With minor modifications, this algorithm can also be used to construct the GLM FLiRTI coefficient path.

4.3. *Model selection.* To fit FLiRTI one must select values for three different tuning parameters, $\lambda$, $\omega$ and the derivative to assume sparsity in, $d$. The choice of $\lambda$ and $d$ in the FLiRTI setting is analogous to the choice of the tuning parameters in a standard smoothing situation. In the smoothing situation one observes $n$ pairs of $(x_i, y_i)$ and chooses $g(t)$ to minimize

$$(19) \qquad \sum_i (y_i - g(x_i))^2 + \lambda \int \{g^{(d)}(t)\}^2 \, dt.$$

In this case the second term, which controls the smoothness of the curve, involves two tuning parameters, namely $\lambda$ and the derivative, $d$. The choice of $d$ in the smoothing situation is completely analogous to the choice of the derivative in FLiRTI. As with FLiRTI, different choices will produce different shapes. A significant majority of the time $d$ is set to 2 resulting in a cubic spline. If the data is used to select $d$ then the most common approach is to choose the values of $\lambda$ and $d$ that produce the lowest cross-validated residual sum of squares.

We adopt the later approach with FLiRTI. In particular we implement FLiRTI using two derivatives, the zeroth and a second derivative, $d$ with $d$ typically chosen from the values $d = 1, 2, 3, 4$. We then compute the cross-validated residual sum of squares for $d = 1, 2, 3, 4$ and a grid of values for $\lambda$ and $\omega$. The final tuning parameters are those corresponding to the lowest cross-validated value. Even though this approach involves three tuning parameters, it is still computationally feasible because there are only a few values of $d$ to test out and, in practice, the results are relatively insensitive to the exact value of $\omega$ so only a few grids points need to be considered. Additional computational savings are produced if one sets $\omega$ to zero. This has the effect of reducing the number of tuning parameters to two by restricting FLiRTI to assume sparsity in only one derivative. We explore this option in the simulation section below and show that this restriction can often still produce good results.

**5. Simulation study.** In this section, we use a comprehensive simulation study to demonstrate four versions of the FLiRTI method, and compare the results with the basis approach using B-spline bases. The first two versions of FLiRTI, "FLiRTI$_L$" and "FLiRTI$_D$," respectively use the lasso and Dantzig methods assuming sparsity in the zeroth and one other derivative. The second two versions, "FLiRTI$_{1L}$" and "FLiRTI$_{1D}$," do not assume sparsity in the zeroth derivative but are otherwise identical to the first two implementations. We consider three cases. The details of the simulation models are as follows.

- Case I: $\beta(t) = 0$ (no signal).
- Case II: $\beta(t)$ is piecewise quadratic with a "flat" region (see Figure 1). Specifically,

$$\beta(t) = \begin{cases} (t - 0.5)^2 - 0.025, & \text{if } 0 \le t < 0.342, \\ 0, & \text{if } 0.342 \le t \le 0.658, \\ -(t - 0.5)^2 + 0.025, & \text{if } 0.658 < t \le 1. \end{cases}$$

- Case III: $\beta(t)$ is a cubic curve (see Figure 3), that is,

$$\beta(t) = t^3 - 1.6t^2 + 0.76t + 1, \qquad 0 \le t \le 1.$$

This is a model where one might not expect FLiRTI to have any advantage over the B-spline method.

In each case, we consider three different types of $X(t)$.

- Polynomials: $X(t) = a_0 + a_1 t + a_2 t^2 + a_3 t^3, 0 \le t \le 1$.
- Fourier: $X(t) = a_0 + a_1 \sin(2\pi t) + a_2 \cos(2\pi t) + a_3 \sin(4\pi t) + a_4 \cos(4\pi t), 0 \le t \le 1$.
- B-splines: $X(t)$ is a linear combination of cubic B-splines, with knots at $1/7, \ldots, 6/7$.

The coefficients in $X(t)$ are generated from the standard normal distribution. The error term $\varepsilon$ in (1) follows a normal distribution $N(0, \sigma^2)$, where $\sigma^2$ is set equal to 1 for the first case and appropriate values for other cases such that each of the signal to noise ratios, $\mathrm{Var}(f(X))/\mathrm{Var}(\varepsilon)$, is equal to 4. We generate $n = 200$ training observations from each of the above models, along with 10,000 test observations.

TABLE 1
*The columns are for different methods. The rows are for different $X(t)$. The numbers outside the parentheses are the average MSEs over 100 repetitions, and the numbers inside the parentheses are the corresponding standard errors*

|  | **B-spline** | **FLiRTI$_L$** | **FLiRTI$_D$** | **FLiRTI$_{1L}$** | **FLiRTI$_{1D}$** |
|---|---|---|---|---|---|
| Case I ($\times 10^{-2}$) | | | | | |
| Polynomial | 2.10 (0.14) | 0.99 (0.14) | 0.38 (0.09) | 1.5 (0.16) | 0.52 (0.11) |
| Fourier | 1.90 (0.18) | 0.53 (0.09) | 0.47 (0.11) | 1.40 (0.15) | 0.50 (0.10) |
| B-spline | 2.40 (0.32) | 0.82 (0.14) | 0.45 (0.11) | 1.40 (0.22) | 0.57 (0.14) |
| Case II ($\times 10^{-5}$) | | | | | |
| Polynomial | 1.20 (0.12) | 0.85 (0.09) | 0.72 (0.09) | 0.92(0.09) | 0.92 (0.08) |
| Fourier | 3.90 (0.32) | 3.40 (0.27) | 3.30 (0.29) | 3.50 (0.30) | 3.60 (0.28) |
| B-spline | 0.52 (0.03) | 0.44 (0.03) | 0.37 (0.03) | 0.43 (0.03) | 0.46 (0.03) |
| Case III ($\times 10^{-2}$) | | | | | |
| Polynomial | 0.96 (0.11) | 0.60 (0.07) | 0.74 (0.08) | 0.57(0.08) | 0.66 (0.08) |
| Fourier | 0.79 (0.11) | 0.43 (0.05) | 0.62 (0.06) | 0.44 (0.05) | 0.46 (0.06) |
| B-spline | 0.080 (0.007) | 0.066 (0.007) | 0.074 (0.008) | 0.063(0.005) | 0.070 (0.008) |

As discussed in Section 4.3, fitting FLiRTI requires choosing three tuning parameters, $\lambda$, $\omega$ and $d$, the second derivative to penalize. For the B-spline method, the tuning parameters include the order of the B-spline and the number of knots (the location of the knots is evenly spaced between 0 and 1). To ensure a fair comparison between the two methods, for each training data set, we generate a separate validation data set also containing 200 observations. The validation set is then used to select tuning parameters that minimize the validation error. Using the selected tuning parameters, we calculate the mean squared error (MSE) on the test set. We repeat this 100 times and compute the average MSEs and their corresponding standard errors. The results are summarized in Table 1.

As we can see, in terms of prediction accuracy, the FLiRTI methods perform consistently better than the B-spline method. Since the FLiRTI$_1$ methods do not search for "flat" regions their results deteriorated somewhat for Cases I and II over standard FLiRTI, correspondingly the results improve slightly in Case III. However, in all cases all four versions of FLiRTI outperform the B-spline method. This is particularly interesting for Case III. Both FLiRTI and the B-spline method can potentially model Case III exactly but only if the correct value for $d$ is chosen in FLiRTI and the correct order in the B-spline. Since these tuning parameters are chosen automatically using a separate validation set neither method has an obvious advantage yet

*The table contains the percentage of truly identified zero region. The rows are for different $X(t)$. The numbers outside the parentheses are the averages over 100 repetitions, and the numbers inside the parentheses are the corresponding standard errors*

|  | FLiRTI$_L$ | FLiRTI$_D$ |
|---|---|---|
| Case I |  |  |
|    Polynomial | 61% (6%) | 65% (6%) |
|    Fourier | 91% (2%) | 71% (6%) |
|    B-spline | 54% (6%) | 74% (6%) |
| Case II |  |  |
|    Polynomial | 70% (6%) | 70% (6%) |
|    Fourier | 72% (5%) | 72% (5%) |
|    B-spline | 58% (5%) | 59% (5%) |

FLiRTI still outperforms the B-spline approach. The lasso and the Dantzig selector implementations of FLiRTI perform similarly, with the Dantzig selector having an edge in the first case, while lasso has a slight advantage in the third case. Finally, for Cases I and II, we also computed the fraction of the zero regions that FLiRTI correctly identified. The results are presented in Table 2.

**6. Canadian weather data.** In this section we demonstrate the FLiRTI approach on a classic functional linear regression data set. The data consisted of one year of daily temperature measurements from each of 35 Canadian weather stations. Figure 4(a) illustrates the curves for 9 randomly selected stations. We also observed the total annual rainfall, on the log scale, at each weather station. The aim was to use the temperature curves to predict annual rainfall at each location. In particular, we were interested in identifying the times of the year that have an effect on rainfall. Previous research suggested that temperatures in the summer months may have little or no relationship to rainfall whereas temperatures at other times do have an effect. Figure 4(b) provides an estimate for $\beta(t)$ achieved using the B-spline basis approach outlined in the previous section. In this case we restricted the values at the start and the end of the year to be equal and chose the dimension, $q = 4$, using cross-validation. The curve suggests a positive relationship between temperature and rainfall in the fall months and a negative relationship in the spring. There also appears to be little relationship during the summer months. However, because of the restricted functional form of the curve, there are only two points, where $\beta(t) = 0$.

The corresponding estimate from the FLiRTI approach, after dividing the yearly data into a grid of 100 equally spaced points and restricting the zeroth and third derivatives, is presented in Figure 4(c) (black line) with the spline estimate in grey. The choice of $\lambda$ and $\omega$ were made using tenfold cross-validation. The FLiRTI estimate also indicates a negative relationship in the spring and a positive relationship in the late fall but no relationship in the summer and winter months. In comparing the B-spline and FLiRTI fits, both are potentially reasonable. The B-spline fit suggests a possible cos/sin relationship, which seems sensible given that the climate pattern is usually seasonal. Alternatively, the FLiRTI fit produces a simple and easily interpretable result. In this example, the FLiRTI estimate seemed to be slightly more accurate with 10 fold cross-validated sum of squared errors of 4.77 vs 5.70 for the B-spline approach.

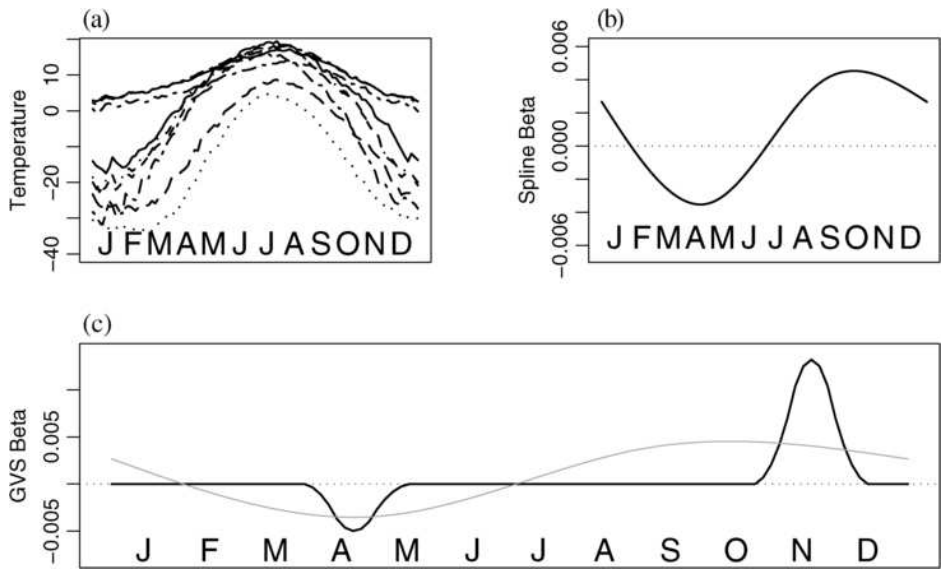

Fig. 4. (a) *Smoothed daily temperature curves for 9 of 35 Canadian weather stations.* (b) *Estimated beta curve using a natural cubic spline.* (c) *Estimated beta curve using FLiRTI approach (black) with cubic spline estimate (grey).*

In addition to estimates for $\beta(t)$, one can also easily generate confidence intervals and tests of significance. We illustrate these ideas in Figure 5. Pointwise confidence intervals on $\beta(t)$ can be produced by bootstrapping the pairs of observations $\{Y_i, X_i(t)\}$, reestimating $\beta(t)$ and then taking the appropriate empirical quantiles from the estimated curves at each time point. Figures 5(a) and (b) illustrate the estimates from restricting the first and third derivatives, respectively, along with the corresponding 95% confidence intervals. In both cases the confidence intervals confirm the statistical significance of the positive relationship in the fall months. The significance of the negative relationship in the spring months is less clear since the upper bound is at zero. However, this is somewhat misleading because approximately 96% of the bootstrap curves did include a dip in the spring but, because the dips occurred at slightly different times, their effect canceled out to some extent. Some form of curve registration may be appropriate but we will not explore that here. Note that the bootstrap estimates also consistently estimate zero relationship during the summer months providing further evidence that there is little effect from temperature in this period. Finally, Figure 5(c) illustrates a permutation test we developed for testing statistical significance of the relationship between temperature and rainfall. The grey line indicates the value of $R^2$ (0.73) for the FLiRTI method applied to the weather data. We then permuted the response variable 500 times and for each permutation computed the new $R^2$. All 500 permuted $R^2$'s were well below 0.73, providing very strong evidence of a true relationship.

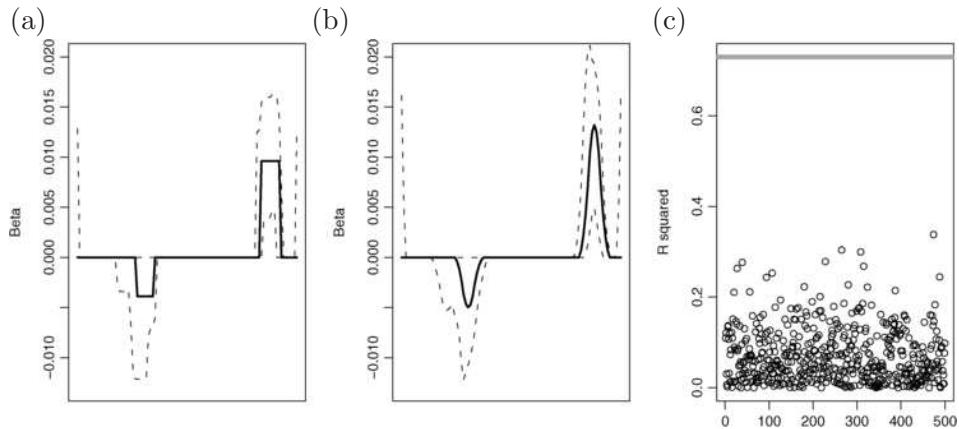

Fig. 5.    (a) *Estimated beta curve from constraining the zeroth and first derivatives.* (b) *Estimated beta curve from constraining the zeroth and third derivatives. The dashed lines represent 95% confidence intervals.* (c) $R^2$ *from permuting the response variable* 500 *times. The grey line represents the observed* $R^2$ *from the true data.*

**7. Discussion.**    The approach presented in this paper takes a departure from the standard regression paradigm, where one generally attempts to minimize an $L_2$ quantity, such as the sum of squared errors, subject to an additional penalty term. Instead we attempt to find the sparsest solution, in terms of various derivatives of $\beta(t)$, subject to the solution providing a reasonable fit to the data. By directly searching for sparse solutions we are able to produce estimates that have far simpler structure than that from traditional methods while still maintaining the flexibility to generate smooth coefficient curves when required/desired. The exact shape of the curve is governed by the choice of derivatives to constrain, which is analogous to the choice of the derivative to penalize in a traditional smoothing spline. The final choice of derivatives can be made either on subjective grounds, such as the tradeoff between interpretability and smoothness, or using an objective criteria, such as the derivative producing the lowest cross validated error. The theoretical bounds derived in Section 3, which show the error rate can grow as slowly as $\sqrt{\log p}$, as well as the empirical results, suggest that one can choose an extremely flexible basis, in terms of a large value for $p$, without sacrificing prediction accuracy.

There has been some previous work along these lines. For example, Tibshirani et al. (2005) uses an $L_1$ lasso-type penalty on both the zeroth and first derivatives of a set of coefficients to produce an estimate which is both exactly zero over some regions and exactly constant over other regions. Valdes–Sosa et al. (2005) also uses a combination of both $L_1$ and $L_2$ penalties on fMRI data. Probably the work closest to ours is a recent approach by Lu and Zhang (2008), called the "functional smooth lasso" (FSL), that was completed at the same time as FLiRTI. The FSL uses a lasso-type approach by placing an $L_1$ penalty on the zeroth derivative and an $L_2$ penalty on the second derivative. This is a nice approach and, as with FLiRTI, produces regions where $\beta(t)$ is exactly zero. However, our approach can be differentiated from these other methods in that we consider derivatives of different order so FLiRTI can generate piecewise constant, linear and quadratic sections. In addition FLiRTI, possesses interesting theoretical properties in terms of the nonasymptotic bounds.

## APPENDIX A: PROOF OF THEOREM 1

We first present definitions of $\delta, \theta$ and $C_{n,p}(t)$. The definitions of $\delta$ and $\theta$ were first introduced in Candes and Tao (2005).

DEFINITION 1.   Let $X$ be an $n$ by $p$ matrix and let $X_T, T \subset \{1, \ldots, p\}$ be the $n$ by $|T|$ submatrix obtained by standardizing the columns of $X$ and extracting those corresponding to the indices in $T$. Then we define $\delta_S^X$ as the smallest quantity such that $(1 - \delta_S^X)\|\mathbf{c}\|_2^2 \leq \|X_T\mathbf{c}\|_2^2 \leq (1 + \delta_S^X)\|\mathbf{c}\|_2^2$ for all subsets $T$ with $|T| \leq S$ and all vectors $\mathbf{c}$ of length $|T|$.

DEFINITION 2.   Let $T$ and $T'$ be two disjoint sets with $T, T' \subset \{1, \ldots, p\}$, $|T| \leq S$ and $|T'| \leq S'$. Then, provided $S + S' \leq p$, we define $\theta_{S,S'}^X$ as the smallest quantity such that $|(X_T\mathbf{c})^T X_{T'}\mathbf{c}'| \leq \theta_{S,S'}^X \|\mathbf{c}\|_2\|\mathbf{c}'\|_2$ for all $T$ and $T'$ and all corresponding vectors $\mathbf{c}$ and $\mathbf{c}'$.

Finally, let $C_{n,p}(t) = \frac{4\alpha_{n,p}(t)}{1 - \delta_{2S}^V - \theta_{S,2S}^V}$, where

$$\alpha_{n,p}(t) = \sqrt{\sum_{j=1}^p \frac{(\mathbf{B}_p(t)^T A_j^{-1})^2}{1/n \sum_{i=1}^n (\int X_i(s)\mathbf{B}_p(s)^T A_j^{-1} ds)^2}}.$$

Next, we state a lemma which is a direct consequence of Theorem 2 in James, Radchenko and Lv (2009) and Theorem 1.1 in Candes and Tao (2007). The lemma is utilized in the proof of Theorem 1.

LEMMA 1. *Let* $\mathbf{Y} = \tilde{X}\tilde{\gamma} + \varepsilon$, *where* $\tilde{X}$ *has norm one columns. Suppose that* $\tilde{\gamma}$ *is an* $S$-*sparse vector with* $\delta_{2S}^X + \theta_{S,2S}^X < 1$. *Let* $\hat{\tilde{\gamma}}$ *be the corresponding solution from the Dantzig selector or the lasso. Then* $\|\hat{\tilde{\gamma}} - \tilde{\gamma}\| \leq \frac{4\lambda\sqrt{S}}{1-\delta_{2S}^X - \theta_{S,2S}^X}$ *provided that* (10) *and* $\max |\tilde{X}^T\varepsilon| \leq \lambda$ *both hold.*

Lemma 1 extends Theorem 1.1 in Candes and Tao (2007), which deals only with the Dantzig selector, to the lasso. Now we provide the proof of Theorem 1. First note that the functional linear regression model given by (6) can be reexpressed as,

$$(20) \qquad \mathbf{Y} = V\gamma + \varepsilon^* = \tilde{V}\tilde{\gamma} + \varepsilon^*,$$

where $\tilde{\gamma} = D_v\gamma$ and $D_v$ is a diagonal matrix consisting of the column norms of $V$. Hence, by Lemma 1, $\|D_v\hat{\gamma} - D_v\gamma\| = \|\hat{\tilde{\gamma}} - \tilde{\gamma}\| \leq \frac{4\lambda\sqrt{S}}{1-\delta_{2S}^V - \theta_{S,2S}^V}$ provided (11) holds.

But $\hat{\beta}(t) = \mathbf{B}_p(t)^T A^{-1}\hat{\gamma} = \mathbf{B}_p(t)^T A^{-1} D_v^{-1}\hat{\tilde{\gamma}}$, while $\beta(t) = \mathbf{B}_p(t)^T\eta + e_p(t) = \mathbf{B}_p(t)^T A^{-1} D_v^{-1}\tilde{\gamma} + e_p(t)$. Then

$$\begin{aligned}
|\hat{\beta}(t) - \beta(t)| &\leq |\hat{\beta}(t) - \mathbf{B}_p(t)^T\eta| + |e_p(t)| \\
&= \|\mathbf{B}_p(t)^T A^{-1} D_v^{-1}(\hat{\tilde{\gamma}} - \tilde{\gamma})\| + |e_p(t)| \\
&\leq \|\mathbf{B}_p(t)^T A^{-1} D_v^{-1}\| \cdot \|\hat{\tilde{\gamma}} - \tilde{\gamma}\| + \omega_p \\
&= \frac{1}{\sqrt{n}}\alpha_{n,p}(t)\|\hat{\tilde{\gamma}} - \tilde{\gamma}\| + \omega_p \\
&\leq \frac{1}{\sqrt{n}}\frac{4\alpha_{n,p}(t)\lambda\sqrt{S_p}}{1-\delta_{2S_p}^V - \theta_{S_p,2S_p}^V} + \omega_p.
\end{aligned}$$

## APPENDIX B: PROOF OF THEOREM 2

Substituting $\lambda = \sigma_1\sqrt{2(1+\phi)\log p} + M\omega\sqrt{n}$ into (12) gives (13). Let $\varepsilon_i' = \int X_i(t)e_p(t)\,dt$. Then to show that (11) holds with the appropriate probability note that

$$\begin{aligned}
|\tilde{V}_j^T\varepsilon^*| = |\tilde{V}_j^T\varepsilon + \tilde{V}_j^T\varepsilon'| &\leq |\tilde{V}_j^T\varepsilon| + |\tilde{V}_j^T\varepsilon'| \\
&= \sigma_1|Z_j| + |\tilde{V}_j^T\varepsilon'| \leq \sigma_1|Z_j| + M\omega\sqrt{n},
\end{aligned}$$

where $Z_j \sim N(0, 1)$. This result follows from the fact that $\tilde{V}_j$ is norm one and, since $\varepsilon_i \sim N(0, \sigma_1)$, it will be the case that $\tilde{V}_j^T \boldsymbol{\varepsilon} \sim N(0, \sigma_1)$. Hence

$$
\begin{aligned}
P\Big(\max_j |\tilde{V}_j^T \boldsymbol{\varepsilon}^*| > \lambda\Big) &= P\Big(\max_j |\tilde{V}_j^T \boldsymbol{\varepsilon}^*| > \sigma_1 \sqrt{2(1+\phi)\log p} + M\omega\sqrt{n}\Big) \\
&\le P\Big(\max_j |Z_j| > \sqrt{2(1+\phi)\log p}\Big) \\
&\le p \frac{1}{\sqrt{2\pi}} \exp\{-(1+\phi)\log p\} / \sqrt{2(1+\phi)\log p} \\
&= (p^\phi \sqrt{4(1+\phi)\pi \log p})^{-1}.
\end{aligned}
$$

The penultimate line follows from the fact that $P(\sup_j |Z_j| > u) \le \frac{p}{u} \frac{1}{\sqrt{2\pi}} \times \exp(-u^2/2)$.

In the case, where $\varepsilon_i^* \sim N(0, \sigma_2^2)$ then substituting $\lambda = \sigma_2 \sqrt{(1+\phi)\log p}$ into (12) gives (14). In this case $\tilde{V}_j^T \boldsymbol{\varepsilon}^* = \sigma_2 Z_j$, where $Z_j \sim N(0, 1)$. Hence

$$
P\Big(\max_j |\tilde{V}_j^T \boldsymbol{\varepsilon}^*| > \sigma_2 \sqrt{(1+\phi)\log p}\Big) = P\Big(\max_j |Z_j| > \sqrt{2(1+\phi)\log p}\Big)
$$

and the result follows in the same manner as above.

## APPENDIX C: THEOREMS 5 AND 6 ASSUMING GAUSSIAN $\varepsilon^*$

The following theorems hold with $\lambda = \sigma_2 \sqrt{2(1+\phi)\log p}$.

THEOREM 5. *Suppose A-1 through A-5 all hold and $\varepsilon_i^* \sim N(0, \sigma_2^2)$. Then, with arbitrarily high probability,*

$$
|\hat{\beta}_n(t) - \beta(t)| \le O(n^{-1/2}) + \frac{H}{p^{*m}}
$$

*and*

$$
\sup_t |\hat{\beta}_n(t) - \beta(t)| \le O(n^{-1/2}) + \frac{H}{p^{*m}}.
$$

Theorem 5 demonstrates that, with the additional assumption that $\varepsilon_i^* \sim N(0, \sigma_2^2)$, asymptotically the approximation error is now bounded by $H/p^{*m}$, which is strictly less than $E_n(t)$. Finally, Theorem 6 allows $p$ to grow with $n$, which removes the bias term, and hence $\hat{\beta}_n(t)$ becomes a consistent estimator.

THEOREM 6. *Suppose we assume A-6 as well as $\varepsilon_i^* \sim N(0, \sigma_2^2)$. Then if we let $p$ grow at the rate of $n^{1/(2m+2b_t)}$, the rate of convergence improves to*

$$
|\hat{\beta}_n(t) - \beta(t)| = O\Big(\frac{\sqrt{\log n}}{n^{1/2(m/(m+b_t))}}\Big)
$$

*or if we let $p$ grow at the rate of $n^{1/(2m+2c)}$, the supremum converges at a rate of*

$$\sup_t |\widehat{\beta}(t) - \beta(t)| = O\left(\frac{\sqrt{\log n}}{n^{1/2(m/(m+c))}}\right).$$

## APPENDIX D: PROOFS OF THEOREMS 3–6

PROOF OF THEOREM 3. By Theorem 2, for $p = p^*$,

$$|\widehat{\beta}(t) - \beta(t)| \leq \frac{1}{\sqrt{n}} C_{n,p^*}(t) \sigma_1 \sqrt{2S_{p^*}(1+\phi)\log p^*} + \omega_{p^*}\{1 + C_{n,p^*}(t) M \sqrt{S_{p^*}}\}$$

with arbitrarily high probability provided $\phi$ is large enough. But, by A-1, $S_{p^*}$ is bounded, and, by A-3 and A-5, $C_{n,p^*}(t)$ is bounded for large $n$. Hence, since $p^*, \sigma_1$ and $\phi$ are fixed, the first term on the right-hand side is $O(n^{-1/2})$. Finally, by A-2, $\omega_{p^*} \leq H/p^{*^m}$ and, by A-1, $S_{p^*} \leq S$ so the second term of the equation is at most $E_n(t)$. The result for $\sup_t |\widehat{\beta}(t) - \beta(t)|$ can be proved in an identical fashion by replacing A-3 by A-4. □

PROOF OF THEOREM 4. By A-1 there exists $S < \infty$ such that $S_p < S$ for all $p$. Hence, by (13), setting $\phi = 0$, with probability converging to one as $p \to \infty$,

$$\begin{aligned}
|\widehat{\beta}(t) - \beta(t)| &\leq \frac{1}{\sqrt{n}} \frac{4\alpha_{n,p_n}(t)}{1 - \delta_{2S}^V - \theta_{S,2S}^V} \sigma_1 \sqrt{2S\log p_n} \\
&\quad + \omega_{p_n}\left\{1 + \frac{4\alpha_{n,p_n}(t)}{1 - \delta_{2S}^V - \theta_{S,2S}^V} M \sqrt{S}\right\} \\
&= \frac{p_n^{b_t}\sqrt{\log p_n}}{\sqrt{n}} \frac{4p_n^{-b_t}\alpha_{n,p_n}(t)}{1 - \delta_{2S}^V - \theta_{S,2S}^V} \sigma_1 \sqrt{2S} \\
&\quad + \frac{\omega_{p_n}p_n^m}{p_n^{m-b_t}}\left\{p_n^{-b_t} + \frac{4p_n^{-b_t}\alpha_{n,p_n}(t)}{1 - \delta_{2S}^V - \theta_{S,2S}^V} M \sqrt{S}\right\} \\
&= \frac{\sqrt{\log n}}{n^{1/2 - b_t/(2m)}} K,
\end{aligned}$$

where $K$ is

$$\begin{aligned}
(21) \quad &\left(\frac{p_n}{n^{1/2m}}\right)^{b_t} \sqrt{\frac{\log p_n}{\log n}} \frac{4p_n^{-b_t}\alpha_{n,p_n}(t)}{1 - \delta_{2S}^V - \theta_{S,2S}^V} \sigma_1 \sqrt{2S} \\
&\quad + \left(\frac{p_n}{n^{1/2m}}\right)^{b_t - m} \frac{\omega_{p_n}p_n^m}{\sqrt{\log n}}\left\{p_n^{-b_t} + \frac{4p_n^{-b_t}\alpha_{n,p_n}(t)}{1 - \delta_{2S}^V - \theta_{S,2S}^V} M \sqrt{S}\right\}.
\end{aligned}$$

But if we let $p_n = O(n^{1/2m})$ then (21) is bounded because $p_n/n^{1/2m}$ and $\log p_n / \log n$ are bounded by construction of $p_n$, $\omega_{p_n} p_n^m$ is bounded by A-2, $p_n^{-b_t} \alpha_{n,p_n}(t)$ is bounded by A-3 and $(1 - \delta_{2S}^V - \theta_{S,2S}^V)^{-1}$ is bounded by A-6. Hence $|\widehat{\beta}_n(t) - \beta(t)| = O(\frac{\sqrt{\log n}}{n^{1/2 - b_t/(2m)}})$. With the addition of A-4 exactly the same argument can be used to prove $\sup_t |\widehat{\beta}(t) - \beta(t)| = O(\frac{\sqrt{\log n}}{n^{1/2 - c/(2m)}})$.   □

PROOF OF THEOREM 5.   By Theorem 2, if we assume $\varepsilon_i^* \sim N(0, \sigma_2^2)$, then, for $p = p^*$,

$$|\widehat{\beta}(t) - \beta(t)| \leq \frac{1}{\sqrt{n}} C_{n,p^*}(t) \sigma_2 \sqrt{2 S_{p^*}(1 + \phi) \log p^*} + \omega_{p^*}$$

with arbitrarily high probability provided $\phi$ is large enough. Then we can show that the first term is $O(n^{-1/2})$ in exactly the same fashion as for the proof of Theorem 3. Also, by A-2, the second term is bounded by $H/p^{*^m}$. □

PROOF OF THEOREM 6.   If $\varepsilon_i^* \sim N(0, \sigma_2^2)$ then, by (14), setting $\phi = 0$, with probability converging to one as $p \to \infty$,

$$|\widehat{\beta}(t) - \beta(t)| \leq \frac{\sqrt{\log n}}{n^{1/2(m/(m+b_t))}} K,$$

where

(22)
$$\begin{aligned} K &= \left( \frac{p_n}{n^{(1/(2(m+b_t)))}} \right)^{b_t} \sqrt{\frac{\log p_n}{\log n}} \frac{4 p_n^{-b_t} \alpha_{n,p_n}(t)}{1 - \delta_{2S}^V - \theta_{S,2S}^V} \sigma_2 \sqrt{2S} \\ &\quad + \left( \frac{p_n}{n^{(1/(2(m+b_t)))}} \right)^{-m} \frac{p_n^m \omega_p}{\sqrt{\log n}}. \end{aligned}$$

Hence if $p_n = O(n^{1/2(m+b)})$ then (22) is bounded, using the same arguments as with (21), so $|\widehat{\beta}_n(t) - \beta(t)| = O(\frac{\sqrt{\log n}}{n^{1/2(m/(m+b_t))}})$. We can prove that $\sup_t |\widehat{\beta}_n(t) - \beta(t)| = O(\frac{\sqrt{\log n}}{n^{1/2(m/(m+c))}})$ in the same way.   □

## APPENDIX E: PROOF OF COROLLARY 1

Here we assume that $A$ is produced using the standard second difference matrix,

(23)
$$A = p^2 \begin{bmatrix} 1/p^2 & 0 & 0 & 0 & \dots & 0 & 0 & 0 \\ -1/p & 1/p & 0 & 0 & \dots & 0 & 0 & 0 \\ 1 & -2 & 1 & 0 & \dots & 0 & 0 & 0 \\ 0 & 1 & -2 & 1 & \dots & 0 & 0 & 0 \\ \vdots & \vdots & \vdots & \vdots & \ddots & \vdots & \vdots & \vdots \\ 0 & 0 & 0 & 0 & \dots & 1 & -2 & 1 \end{bmatrix}.$$

Throughout this proof let $\eta_k = \beta(k/p)$. First we show A-1 holds. Suppose that $\beta''(t) = 0$ for all $t$ in $R_{k-2}, R_{k-1}$ and $R_k$ then there exist $b_0$ and $b_1$ such that $\beta(t) = b_0 + b_1 t$ over this region. Hence, for $k \geq 2$,

$$
\begin{aligned}
\gamma_k &= p^2(\eta_{k-2} - 2\eta_{k-1} + \eta_k) \\
&= p^2(\beta((k-2)/p) - 2\beta((k-1)/p) + \beta(k/p)) \\
&= p^2\left(b_0 + b_1\frac{k-2}{p} - 2b_0 - 2b_1\frac{k-1}{p} + b_0 + b_1\frac{k}{p}\right) = 0.
\end{aligned}
$$

But note that if $\beta''(t) \neq 0$ at no more than $S$ values of $t$ then there are at most $3S$ triples such that $\beta''(t) \neq 0$ for some $t$ in $R_{k-2}, R_{k-1}$ and $R_k$. Hence there can be no more than $3S + 2$ $\gamma_k$'s that are not equal to zero (where the two comes from $\gamma_1$ and $\gamma_2$).

Next we show A-2 holds. For any $t \in R_k$, $\mathbf{B}(t)^T \boldsymbol{\eta} = \eta_k$. But since $|\beta'(t)| < G$ for some $G < \infty$ and $R_k$ is of length $1/p$ it must be the case that $\sup_{t \in R_k} \beta(t) - \inf_{t \in R_k} \beta(t) \leq G/p$. Let $\eta_k$ be any value between $\sup_{t \in R_k} \beta(t)$ and $\inf_{t \in R_k} \beta(t)$ for $k = 1, \ldots, p$. Then

$$
\omega_p = \sup_t |\beta(t) - \mathbf{B}(t)^T \hat{\boldsymbol{\eta}}| = \max_k \sup_{t \in R_k} |\beta(t) - \eta_k|
$$

$$
\leq \max_k \left\{ \sup_{t \in R_k} \beta(t) - \inf_{t \in R_k} \beta(t) \right\} \leq G/p
$$

so A-2 holds with $m = 1$.

Now, we show A-3 holds. For $t \in R_k$ let $L_{nj}(t) = \frac{1}{n}\sum_{i=1}^n(\frac{1}{p}\sum_{l=1}^p X_{il}\frac{A_{lj}^{-1}}{A_{kj}^{-1}})^2$, where $A_{lk}^{-1}$ is the $l, k$th element of $A^{-1}$ and $X_{il}$ is the average of $X_i(t)$ in $R_l$, that is, $p\int_{R_l} X_i(s)\,ds$. Then $\alpha_{n,p}(t) = \sqrt{\sum_{j=1}^p L_{nj}(t)^{-1}}$. Since $X_i(t)$ is bounded above zero $L_{nj}(t) \geq W_\varepsilon^2(\frac{1}{p}\sum_{l=1}^p\frac{A_{lj}^{-1}}{A_{kj}^{-1}})^2$ for some $W_\varepsilon > 0$. It is easily verified that

(24) $$
A^{-1} = \frac{1}{p^2}\begin{bmatrix} p^2 & 0 & 0 & 0 & \ldots & 0 \\ p^2 & p & 0 & 0 & \ldots & 0 \\ p^2 & 2p & 1 & 0 & \ldots & 0 \\ p^2 & 3p & 2 & 1 & \ldots & 0 \\ \vdots & \vdots & \vdots & \vdots & \ddots & 0 \\ p^2 & p(p-1) & p-2 & p-3 & \ldots & 1 \end{bmatrix}
$$

and, hence,

$$
\frac{A_{lj}^{-1}}{A_{kj}^{-1}} = \begin{cases} \infty, & k < j,\, l \geq j, \\ 0, & k \geq j,\, l < j, \\ \dfrac{l-j+1}{k-j+1}, & k \geq j,\, l \geq j \end{cases}
$$

except for $j = 1$ in which case the ratio equals one for all $l$ and $k$. For $t = 0$ then $t \in R_1$ (i.e., $k = 1$) and, hence, $L_{n1}(0) \geq W_\varepsilon^2$ while $L_{nj}(0) = \infty$ for all $j > 1$. Therefore, $\alpha_{n,p}(0) \leq W_\varepsilon^{-1}$ for all $p$. Hence, A-3 holds with $b_0 = 0$. Alternatively, for $0 < t < 1$, then $k = \lfloor pt \rfloor$ and $L_{nj}(t) = \infty$ for $j > k$. Hence for $j \leq k$,

$$L_{nj}(t) \geq W_\varepsilon^2 \frac{1}{p^2} \left( \sum_{l=j}^{p} \frac{l-j+1}{k-j+1} \right)^2$$

$$= W_\varepsilon^2 \frac{1}{p^2} \left( \frac{(p-j+1)(p-j+2)}{2(k-j+1)} \right)^2$$

$$\geq W_\varepsilon^2 \frac{1}{p^2} p^4 \frac{(1-t)^4}{4} \left( \frac{1}{k-j+1} \right)^2 \qquad \text{since } j \leq pt.$$

Therefore

$$\alpha_{n,p}(t) \leq \frac{2}{(1-t)^2 W_\varepsilon p} \sqrt{\sum_{j=1}^{k} (k-j+1)^2}$$

(25)
$$= \frac{2}{(1-t)^2 W_\varepsilon p} \sqrt{\frac{k}{6}(k+1)(2k+1)}$$

$$\leq p^{1/2} \frac{2}{(1-t)^2 W_\varepsilon} t^{3/2}$$

since $k = \lfloor pt \rfloor$. Hence A-3 holds with $b_t = 1/2$ for $0 < t < 1$.

Finally, note that (25) holds for any $t < 1$ and is increasing in $t$ so

$$\sup_{0 < t < 1-a} \alpha_{n,p}(t) \leq p^{1/2} \frac{2}{a W_\varepsilon}$$

for any $a > 0$ and hence A-4 holds with $c = 1/2$.

**Acknowledgments.** We would like to thank the Associate Editor and referees for many helpful suggestions that improved the paper.

G. M. James
Marshall School of Business
University of Southern California
E-mail: gareth@usc.edu

J. Wang
Department of Statistics
University of Michigan
E-mail: jingjw@umich.edu

J. Zhu
Department of Statistics
University of Michigan
E-mail: jizhu@umich.edu